

\baselineskip=14pt
\parskip=10pt

\font\eightrm=cmr8 

\magnification=\magstephalf

\def\1{{\overline{1}}}
\def\2{{\overline{2}}}
\parindent=0pt
\overfullrule=0in

\def\frac#1#2{{#1 \over #2}}
\centerline
{
\bf Searching for Ap\'ery-Style Miracles 
}
\centerline
{\bf [Using, Inter-Alia, the Amazing Almkvist-Zeilberger Algorithm]}
\bigskip
\centerline
{\it Shalosh B. EKHAD and Doron ZEILBERGER}
\bigskip

\qquad \qquad \qquad \qquad {\it Dedicated to Gert Almkvist (b. April 17, 1934), on becoming an octogenerian }

{\bf VERY IMPORTANT}

Like in all the joint articles  of the authors, {\bf this} article is {\bf not} the main point. It may be viewed
as a {\it user's manual} for the much more important {\bf Maple package}, {\tt NesApery}
(BTW, {\it nes} means `miracle' in Hebrew).

It may be obtained directly from the following url:

{\tt http://www.math.rutgers.edu/\~{}zeilberg/tokhniot/NesApery} \quad ,

or via the {\it almost as important} {\bf front} of this article:

{\tt http://www.math.rutgers.edu/\~{}zeilberg/mamarim/mamarimhtml/apery.html} \quad ,

that contains links to {\bf twenty} input and output files, that readers are welcome to use
as templates in order to carry out their own investigations.

{\bf How the Pythagoreans actually proved that the square-root of two is irrational}

The first {\it crisis} in the history of mathematics was when the Pythagoreans, who {\bf correctly}
believed that {\it all} quantities are  ratios of integers, discovered an alleged
`counterexample', namely that the length of the  diagonal of the unit square, that nowadays
is denoted by $\sqrt{2}$, is {\bf not} such a ratio, in other words, is what we call today an
{\it irrational number}. They were upset for nothing! (and should have left poor Hipassus alone).
What they {\it did} prove was that $\sqrt{2}$ does {\bf not} exist (from the philosophical-ontological
viewpoint [of course it exists as a {\it symbol}]), since {\bf by definition}, only integers
and their ratios really exist, and all the rest is the {\it fictional} (and often boring)
creation of humankind.

Anyway, the proof that the ancient Greeks had was {\bf not} the standard textbook proof with which
we torture our students today.

The usual textbook proof uses:
$$
a^2=2b^2 \quad implies  \quad a'^2=2b'^2 \quad ,
$$
where $a'=a/2$ and $b'=b/2$, and since $a'+b'<a+b$ and $a=1,b=0$ is {\bf not} a solution, {\bf qed}.

On the other hand, the original Greek proof, rendered algebraically, is based on
$$
a^2=2b^2 \quad implies \quad  a'^2=2b'^2 \quad ,
$$
where {\bf now} $a'=2b-a$ and $b'=a-b$. Once again $a'+b'=b<a+b$ and we have a {\it reductio} proof.

{\bf How the Pythagoreans should have proved that the square-root of two is irrational}

They should have gone to {\tt Maple} (or Sage, or Mathematica, etc.) and typed

{\tt convert(sqrt(2),confrac); } \quad , 

and would have gotten the output

{\tt [1, 2, 2, 2, 2, 2, 2, 2, 2, 2, 2]} ,

that would have lead them to make the conjecture that
$$
\sqrt{2}= [1,2^{\infty}] \quad,
$$
i.e.
$$
\sqrt{2}=1+1/(2+1/(2+1/(2+1/(2+ \dots \quad, \quad ad \quad infinitum \, .
$$
Since $\sqrt{2}$ is (conjecturally for now) an {\it infinite} simple continued fraction, it
follows immediately that it is {\bf irrational}. Indeed, thanks to their buddy Euclid, who was
going to be  born a couple of centuries later, who made up the {\it Euclidean algorithm}, a rational number is
always expressible as a {\it finite} continued fraction.

Anyway, in order to (rigorously) prove the conjectured form for the continued fraction for $\sqrt{2}$, they
should have first defined
$$
y\, := \, 2+1/(2+1/(2+1/(2+ \dots \quad,
$$
and using {\bf self-similarity}, deduce the equation
$$
y=2+1/y \quad,
$$
that is equivalent to the {\it quadratic equation}
$$
y^2-2y-1=0 \quad .
$$
Then they could have asked their Babylonian neighbors, across the Mediterranean, who were better in algebra, to 
solve this quadratic equation, getting
$$
y=\frac{2 +\sqrt{8}}{2}= 1+ \sqrt{2} \quad .
$$
(Of course $y$ is positive). Then going back to 
$$
x:=1+1/(2+1/(2+1/(2+1/(2+ \dots \quad,
$$
we get
$$
x=1+\frac{1}{y}= 1+\frac{1}{\sqrt{2}+1}= 1+\frac{\sqrt{2}-1}{\sqrt{2}^2-1^2}= 1+(\sqrt{2}-1)= \sqrt{2} \quad .
$$
Yea!, the infinite continued fraction $[1,2^{\infty}]$ indeed equals $\sqrt{2}$ and hence the latter is irrational.

{\bf Why is this proof much better?}

Let's consider, instead of $\sqrt{2}$, $1+\sqrt{2}$. Of course if the latter is irrational, so is the former.

Let's truncate the continued fraction after $n$ $2$'s,  and call it $\frac{p_n}{q_n}$, i.e.
$$
\frac{p_n}{q_n} := [2,2,2, \dots, 2] \quad, \quad (repeated \quad n \quad times) \quad .
$$
Then 
$$
\frac{p_n}{q_n} = 2+\frac{1}{\frac{p_{n-1}}{q_{n-1}}} = 2+\frac{q_{n-1}}{p_{n-1}} = \frac{2p_{n-1}+q_{n-1}}{p_{n-1}} \quad .
$$
Equating numerator and denominator we have
$$
p_n = 2p_{n-1} + q_{n-1} \quad , \quad q_n = p_{n-1} \quad .
$$
Hence
$$
p_n=2p_{n-1}+p_{n-2} \quad ,
$$
and since $q_n=p_{n-1}$ also
$$
q_n=2q_{n-1}+q_{n-2} \quad  .
$$
It turns out that the numerator and denominator sequences of the $n$-th convergent to $1+\sqrt{2}$ {\bf both}
satisfy the  recurrence
$$
x_n=2x_{n-1}+x_{n-2} \quad  ,
$$
but with {\bf different} {\it initial conditions} (of course!, or else their ratio would be $1$) :
$$
p_0=1 \quad, \quad p_1= 2 \quad ;
$$
$$
q_0=0 \quad, \quad q_1= 1 \quad .
$$
So we can {\bf forget} about continued fractions, and just state the following proposition.

{\bf Proposition.} Let $p_n$, $q_n$ be the sequences that are solutions of the linear recurrence
$$
x_n=2x_{n-1}+ x_{n-2} \quad,
$$
with initial conditions $x_0=1,x_1=2$; $x_0=0,x_1=1$, respectively. Then 
$$
\alpha :=\lim_{n \rightarrow \infty} \,\, \frac{p_n}{q_n} \quad ,
$$
exists and is an irrational number.

Of course this is all very classical stuff (see e.g. [N], ch. 5), adapted to fit the more general scenario that
is coming up soon, but still as a motivating example, let's redo the proof, with our new twist
(not using continued fractions explicitly).

Until further notice, $\alpha$ stands for $1+\sqrt{2}$.

{\bf How close is $\frac{p_n}{q_n}$ to $\sqrt{2}+1$?}

Recall that
$$
p_n=2p_{n-1}+p_{n-2} \quad , \quad q_n=2q_{n-1}+q_{n-2} \quad .
$$
Multiplying the first equation by $q_{n-1}$ and the second equation by $p_{n-1}$ we get:
$$
p_n q_{n-1} =2p_{n-1} q_{n-1} +p_{n-2} q_{n-1} \quad , \quad p_{n-1} q_n =2 p_{n-1} q_{n-1}+ p_{n-1} q_{n-2} \quad .
$$
Subtracting the first from the second, we get
$$
p_n q_{n-1} - p_{n-1} q_{n} = (-1) (p_{n-1} q_{n-2} - p_{n-2} q_{n-1}) \quad .
$$
It would be a good idea to give the sequence $p_n q_{n-1} - p_{n-1} q_{n}$ a name, let's call it
$c_n$, to wit:
$$
c_n :=p_n q_{n-1} - p_{n-1} q_{n} \quad .
$$
We have just seen that $c_n$ satisfies the simple {\bf first-order recurrence} (with constant coefficients)
$$
c_n=-c_{n-1}  \quad .
$$
Since $c_1=p_1q_0 - p_0 q_1 = 2 \cdot 0 - 1 \cdot 1 = -1$, we can even {\bf solve} for $c_n$ explicitly and get
the beautiful {\bf explicit} expression
$$
c_n = (-1)^n \quad .
$$

Hence
$$
\frac{p_n}{q_n} - \frac{p_{n-1}}{q_{n-1}}= \frac{c_n}{q_n q_{n-1}} = \frac{(-1)^n}{q_n q_{n-1}}  \quad .
$$
Hence
$$
\frac{p_n}{q_n}= 2 + \sum_{i=2}^{n}  \frac{(-1)^n}{q_n q_{n-1}} \quad .
$$
Since $q_n=O(\alpha^n)$ (why?), it follows that $\lim_{n \rightarrow \infty } \frac{p_n}{q_n}$ {\it exists},
and equals $\alpha=1 + \sqrt{2}$ (why?). Furthermore, the `error', $\alpha -\frac{p_n}{q_n}$, can be easily
estimated. We have
$$
| \, \alpha -\frac{p_n}{q_n} \, | = | \, \sum_{i=n+1}^{\infty}   \frac{(-1)^i}{q_i q_{i-1}} \, | \, \leq \frac{C}{q_n^2} \quad \quad,
$$
for some easily computable constant, $C$, that in this case  may be taken to be $1$.

In particular, we may deduce that $\alpha=1+\sqrt{2}$ is {\it irrational}. Let's see why.

Let's state the easy but beautiful lemma (see, e.g., Ivan Niven's classic [N], Theorem 4.3, and van der Poorten's
expository masterpiece [vdP], p. 196).

{\bf Irrationality Criterion (and measure)}: Given a constant $\alpha$,
if you can exhibit an {\it infinite} sequence of {\it distinct} rational numbers $\frac{p_n}{q_n}$ such that
for {\it some} constant $\delta>1$ and some positive constant $C$, we have
$$
|\alpha - \frac{p_n}{q_n}| \leq \frac{C}{q_n^\delta} \quad,
$$
then $\alpha$ is definitely irrational. Furthermore, if the $q_n$ are of {\it exponential growth}, then
$\alpha$ has so-called {\it irrationality measure} $\leq$ $\delta/(\delta-1)$.

It follows that $1+\sqrt{2}$ (and hence $\sqrt{2}$) is irrational with irrationality measure $2$.

From now on $\alpha$ is no longer $1+ \sqrt{2}$.

Once you know, {\it for sure},  (by {\it whatever} means) that a certain number, let's call it $\alpha$
(no longer (necessarily) $1+ \sqrt{2}$) is irrational, then you can do much better.
Dirichlet famously proved
(see [N], Theorem 4.1), using the (existential!) {\it pigeon principle}, that if $\alpha$ is irrational, then one can always come
up with infinitely many rational numbers $\frac{h}{k}$ such that
$$
|\alpha - \frac{h}{k}| \leq \frac{C}{k^2} \quad \quad ,
$$
where one can take $C=1$. 

{\eightrm
[ In fact (see [N], Th. 5.1), one can even take $C=\frac{1}{\sqrt{5}}$, 
but ([N], Th. 5.2), for $\alpha=\frac{1+\sqrt{5}}{2}$, that $C$ is the smallest possible). ]}

Sometimes you can do even better. A real number $\alpha$ is {\it approximable to order $n$}
if the inequality
$$
| \alpha - \frac{h}{k}| < \frac{C}{k^n} \quad ,
$$
has infinitely many rational solutions $\frac{h}{k}$ with $k>1$ and $gcd(h,k)=1$.

It is easy to see (e.g. Theorem 7.7 of [N], [vdP], section 2) that a rational number is
approximable to order $1$ but {\bf not} higher. So if we have a candidate, $\alpha$,
that we want to prove is irrational, all we have to do is come up with 
(i.e. construct!)  {\it any}  infinite sequence of distinct rational numbers, $\frac{p_n}{q_n}$, such that
we can show that for {\bf some} $\delta>1$ (even, say, $\delta=1+10^{-10^{1000}}$)
$|\alpha - \frac{p_n}{q_n}| \leq \frac{C}{q_n^\delta}$, for {\it some} positive constant $C$.

{\eightrm
Liouville famously proved ([N], Th. 7.8) that an algebraic number of order $n$ is
{\bf not} approximable to order $n+1$ or higher, that lead him to come up
with the first example of {\it transcendental numbers}. The Dyson-Siegel-Thue-Roth
theorem tells you that, in fact, the irrationality measure of any algebraic number is $2$,
but, the proof, just like Dirichlet's, is non-effective, and it is still
a challenge to come-up with non-trivial {\it effective irrationality measures}, even for algebraic numbers.
For example, very smart people (e.g. [Ba], [Be]) had to work very hard in order to get an {\it effective} irrationality measure for
$\root 3 \of 2$ that is strictly less than $3$.}

Once we have the {\it assurance} that our number is indeed irrational, of course, we can do much better by `cheating'. Find
good floating-point approximations, convert to a continued-fraction, and take the
sequence of {\bf convergents}, that would give you the {\bf  best possible} rational approximations.
But before we have the assurance that our candidate constant is indeed irrational, we
have no {\it guarantee} that we can keep going. In other words, the continued fraction
may (at least in principle) be finite, and then we can't go for ever.

{\it Notation Alert}: So far, bowing to tradition, we used {\it subscripts}, $p_n,q_n$ to denote sequences,
being discrete liberationists, we will now revert to functional format, $p(n),q(n),$ etc.

{\bf Man Muss Immer Umkehren (One Must Always Invert)}

It is very possible that someone in a different galaxy, who has never heard of the
process of taking the square-root (but does know how to add and divide) could have
come up with $1+\sqrt{2}$ as the limit of the sequence $\frac{p(n)}{q(n)}$ 
described above.

Let's {\bf generalize}!

Let's {\bf start} with {\it any} homogeneous linear recurrence equation
$$
\sum_{i=0}^{L} a_i(n) x(n+i)=0 \quad ,
\eqno{(GeneralRecurrence)}
$$
where {\bf now} the {\bf coefficients} are allowed to be {\bf polynomials} in $n$, rather than
mere constants (as in the motivating example, where the recurrence was $x(n+2)-2x(n+1)-x(n)=0$).

Define two sequences of rational numbers (not necessarily integers) that satisfy that recurrence
with two different sets of  {\it initial conditions}
$$
p(0)=p_0 \quad , \quad \dots \quad , \quad   p(L-1)=p_{L-1} \quad , \quad 
$$
$$
q(0)=q_0 \quad , \quad \dots \quad , \quad  q(L-1)=q_{L-1} \quad , \quad 
$$
for {\it some} rational numbers (usually integers) $p_0, \dots, p_{L-1} ; q_0, \dots , q_{L-1}$, and {\bf define}
$$
\alpha := \lim_{n \rightarrow \infty} \frac{p(n)}{q(n)} \quad .
$$
Then use the recurrence to crank out the first $300$ (or whatever) terms of $p(n)$ and $q(n)$ and
get an estimate for $\alpha$, say $\frac{p(300)}{q(300)}$.

Then using {\it identification} programs, like LLL, or PSLQ, (e.g. in Maple, the command {\tt identify}),
we ask the computer to {\it conjecture} an alternative description of $\alpha$. If in luck, it will turn out to be
one of the famous constants (and if really lucky, one of the constants for which it is still
open to find a proof of irrationality, e.g. Euler's constant, $\gamma$, or Catalan's constant $C$,
or $e+\pi$, or $e\pi$).

Then you ask your beloved (numeric) computer to estimate a $\delta$ such that
$$
|\alpha - \frac{p(n)}{q(n)}| \leq \frac{C}{q'(n)^{\delta}} \quad ,
$$
where $q'(n)$ is the denominator of the reduced form of $\frac{p(n)}{q(n)}$ (i.e. you write $\frac{p(n)}{q(n)}$ in reduced form $\frac{p'(n)}{q'(n)}$
(with $p'(n), q'(n)$ {\bf integers} and $gcd(p'(n),q'(n))=1$)).

If really {\it in luck}, and $\delta$ would turn out to be larger than $1$, than you are potentially famous.
You have discovered (empirically for now) an {\bf Ap\'ery miracle}.

Having done this preliminary {\it reconnaissance}, we will soon see how to turn this into
a fully rigorous proof, but before that let's describe the three kinds of {\it Ap\'ery miracles}.

{\bf The three kinds of Ap\'ery miracles}

A {\it minor Ap\'ery miracle} is when
$$
\alpha:=\lim_{n \rightarrow \infty} \frac{p(n)}{q(n)} \quad
$$
equals a well-known constant, expressible in terms of $\pi, e, \gamma$, logarithms of 
explicit algebraic numbers, etc.

A {\it major Ap\'ery miracle} is when
$$
\alpha:=\lim_{n \rightarrow \infty} \frac{p(n)}{q(n)} \quad
$$
equals a well-known constant, expressible in terms of $\pi, e$, logarithms of algebraic numbers etc.
{\it and} it appears that $\delta>1$, thereby (tentatively) proving irrationality.

A {\it super Ap\'ery miracle} is when, in addition, the constant, $\alpha$, has not yet been proven to
be irrational, making you (and your computer) (potentially) famous.

Of course, in order to be {\it really} famous, you would have to prove it all rigorously. 
If you are having trouble, you can always try and collaborate with
an expert number theorist.
But in many cases, it is possible, as we will see below, 
by following the same outline as in the simple motivating example above
(where we reproved, via recurrences, that $1+\sqrt{2}$ is irrational),
essentially following   Ap\'ery's proof (as rendered in [vdP]).

{\bf Blessed are the meek}

Often the $\delta$ (both empirical and rigorous) would turn out to be less than $1$.
While it is very disappointing from the {\it number theory} viewpoint, it is
still of some interest  from the {\it numerical analysis}  viewpoint, since
we found a quickly-converging sequence of rational numbers to our constant $\alpha$.

{\bf Operator Notation}

Let $N$ be the {\it shift operator} 
$$
N x(n):=x(n+1) \quad,
$$
then the above, general, recurrence
$$
\sum_{i=0}^{L} a_i(n) x(n+i)=0 \quad,
\eqno{(GeneralRecurrence)}
$$
can be written, in {\it operator notation}, as
$$
\left ( \sum_{i=0}^{L} a_i(n) N^i \right ) \, x(n) \, =0 \quad .
\eqno{(GeneralRecurrence')}
$$
Writing
$$
ope:= \sum_{i=0}^{L} a_i(n) N^i \quad,
$$
let's denote the above-mentioned $\alpha$ and $\delta$ as follows.
$$
\alpha:=RA(ope, [p_0, \dots, p_{L-1}], [q_0, \dots, q_{L-1}]) \quad ,
$$
$$
\delta:=RAdel(ope, [p_0, \dots, p_{L-1}], [q_0, \dots, q_{L-1}]) \quad .
$$

For example,
$$
RA( \, N^2-2N-1 \, , \, [1,2] \, ,\, [0,1]) \, =  \, 1+ \sqrt{2} \quad ,
$$
$$
RAdel(N^2-2N-1 \, , \, [1,2] \, , \, [0,1])= \, 2 \quad .
$$

{\bf The Original Ap\'ery Super Miracle}

$$
RA(\,  \left( n+2 \right)^{3}{N}^{2}-  \left( 2\,n + 3\right) \left( 17\,{n}^{2}+51\,n+39 \right) N  + \left( n+1 \right) ^{3} \, ,\,
[0,6],[1,5])= \zeta(3) \quad , \quad  \delta=1.080529431\dots .
$$
(See [vdP]) [Stand by for an automated (rigorous) proof of this, using Ap\'ery's method].

Indeed, in {\it hindsight}, it {\it could} have been discovered, empirically, by running in {\tt NesApery},
the command

{\tt DelSeqRec(ope,n,N,Ini1,Ini2,K)} \quad ,
for 
$$
ope= \left( n+2 \right)^{3}{N}^{2}-  \left( 2\,n + 3 \right) \left( 17\,{n}^{2}+51\,n+39 \right) N  + \left( n+1 \right) ^{3}
$$
$$
Ini1=[0,6] \quad , \quad Ini2=[1,5] \quad , \quad  K=200 \,\, (say).
$$

This yields that the constant in question is most probably $\zeta(3)$, and that the empirical $\delta$ seems
to be around $1.08$. 

Of course the Ap\'ery $\zeta(3)$ linear recurrence operator is a good one. Most operators, with random initial conditions,
would lead to obscure constants and $\delta$ less than $1$, but at least in principle, (and possibly
also in practice), it would have been possible to come up with that operator by an {\it exhaustive} search.

\vfill\eject

{\bf The Two Other Ap\'ery Miracles}

$$
ope=  \left( n+2 \right) ^{2}{N}^{2} - \left ( 11\,{n}^{2} + 33\,n + 25 \right ) N  - \left ( n+1 \right) ^{2}  \quad , \quad
$$
$$
Ini1=[0,5] \quad, \quad Ini2=[1,3] \quad ,
$$
yielding
$\alpha = \pi^2/6$ (alias $\zeta(2)$) with  $\delta=1.092159255$.
Also:
$$
ope=\left( n+2 \right) {N}^{2}- \left( 9+6\,n \right) N +(n+1) \quad ,
$$
$$
Ini1=[0,2] \quad, \quad Ini2=[1,3] \quad ,
$$
yielding $\alpha=\ln 2$ with $\delta=1.276082872$, (see [vdP]). 
These are {\it major} (but not super) miracles. The constants ($\pi^2$ and $\ln 2$)
are famous, sure enough, and the $\delta$ is larger than $1$, but
since both $\pi^2$ and $\ln 2$ were already proven (for a very long time!) to be irrational, it is not
as exciting.

{\bf Systematic Empirical Search}

Inspired by the $\ln 2$ miracle, procedure

{\tt NisimRec1emp(n,N,K,A,B)}

inputs symbols {\tt n} and {\tt N}, and positive integers {\tt K,A,B}, and outputs {\it all operators} of the form
$$
(n+2)N^2-(a\,n \, + \, b)N+(n+1) \quad ,
$$
with $A \leq a,b \leq B$,
such that, by taking the first $K$ terms, the empirical $\delta$,
$$
RAdel(ope,[0,1],[1,1])
$$
appears to be larger than $1$. With the choice $A=1,B=20$ (i.e. $400$ candidates), the first-named author found
$121$ of them. To see them all, go to:

{\tt http://www.math.rutgers.edu/\~{}zeilberg/tokhniot/oNesApery19} \quad .

The ``champion'' (with the largest empirical $\delta$) happened to be
$$
(n+2)N^2 - (10n+15)N+ n+1 \quad ,
$$ 
with empirical $\delta=1.440802982$ (with $K=100$).

{\bf Alas, the Haystack is Too Big}

Once a genie supplied us with a potential Ap\'ery miracle, it is easy enough to verify,
empirically, that it is indeed a miracle.
But how to come up with them in the first place?

While a pure brute-force search for `cousins' of $\ln 2$ is still feasible, 
such a search for the operators that proved the irrationality of $\zeta(2)$ and $\zeta(3)$
would take too long. But anyone who has read van der Poorten's lovely paper [vdP]
would know that the recurrence operators (or equations) cited above,
that lead to the proofs of $\zeta(3), \zeta(2)$ and $\ln 2$, respectively,
are those annihilating (resp. satisfied by) the {\it hypergeometric sums}
$$
\sum_{k=0}^{n} {{n} \choose {k}}^2 {{n+k} \choose {k}}^2 \quad , \quad
\sum_{k=0}^{n} {{n} \choose {k}}^2 {{n+k} \choose {k}} \quad , \quad
\sum_{k=0}^{n} {{n} \choose {k}} {{n+k} \choose {k}} \quad .
$$
These immediately supply {\it explicit} expressions for the denominators $q(n)$ of the rational approximations $\frac{p(n)}{q(n)}$
for the desired constants. Thanks to the famous {\bf Zeilberger algorithm} ([Ze]) these recurrences could be
derived, in a few seconds, complete with fully rigorous proofs. Then one can use 
a similar proof to the easy proof of the irrationality of $1+ \sqrt{2}$ given above, that is essentially
Ap\'ery's method
(as spelled-out by Don Zagier and Henri Cohen, and so charmingly described  by Alf van der Poorten [vdP]).
This can all be automated, and done by the computer, as we will see below.

This is implemented in procedure {\tt RAsum} in the Maple package {\tt NesApery}. 
A verbose version is {\tt RAsumV}.

For a computerized redux, {\bf in three seconds}, of the above-mentioned Ap\'ery proofs for $\ln 2$, $\zeta(2)$, and $\zeta(3)$, go to:

{\tt http://www.math.rutgers.edu/\~{}zeilberg/tokhniot/oNesApery16} \quad .

But what about other binomial coefficients sums? We can have our computer do the
Ap\'ery treatment to {\bf any} binomial coefficients sum that would describe the
denominators, $q(n)$. Then define the sequence of numerators, $p(n)$, as the one satisfying the
same recurrence, but with initial conditions, say $p(0)=0$, $p(1)=1$ (in the case of
a second order recurrence), or $p(0)=0$, $p(1)=0$, $p(2)=1$ (in the case of third order recurrence), etc.
The only current `gap' in our automated proof, is the statement that

$p(n) lcm(1,...,n)$ (in the case of $\ln 2$) \quad ,  \quad
$p(n) lcm(1,...,n)^2$ (in the case of $\zeta(2)$) \quad ,  \quad
$p(n) lcm(1,...,n)^3$ (in the case of $\zeta(3)$) \quad ,

are integers, but we are sure that these would be easily filled by a sufficiently interested human or machine.

Let's describe procedure {\tt RAsum(F,k,n,N,K)} in some detail, since it is very important.

{\bf Procedure {\tt RAsum(F,k,n,N,K)}}

{\bf Input of {\tt RAsum}}

$\bullet$ A binomial coefficients summand, {\tt F}, 
phrased in terms of the symbols {\tt k} and {\tt n} .

$\bullet$ Symbols {\tt k} and {\tt n} used in {\tt F}  .

$\bullet$ A symbol {\tt N} for the shift operator in  {\tt n} .

$\bullet$ A positive integer {\tt K} (for the upper bound for the order of the recurrence operator
annihilating $c(n):=p(n)q(n-1)-p(n-1)q(n)$, just for the sake of efficiency, one can always make it bigger)

{\bf Output  of {\tt RAsum}}

$\bullet$ The operator, let's call it $ope(n,N)$, annihilating the sequence  defined by the sum, i.e. the denominator sequence $q(n)$, using the
Zeilberger algorithm .

$\bullet$ The operator,  $ope_1(n,N)$,
annihilating  $c(n):=p(n)q(n-1)-p(n-1)q(n)$, it always exists, by general theory, hence
it is fully justifiable to search for it by `guessing' (using undetermined coefficients), for the
famous Ap\'ery cases that operator is always first-order, and $c(n)$ is extremely simple.

$\bullet$ A four-tuple of integers $[G,R,P,D_0]$, with the property that, defining
$$
L(n):=lcm(\, 1, 2, \dots, n) \quad ,
$$
that 
$$
 L(Gn)^R \, P^n \, D_0 \,  \cdot \, p(n)
$$
are always integers. This is the only not-yet-fully-rigorous part, but if needed should be (presumably)
easy to fill, in the original Ap\'ery miracles narrated in [vdP] it was done by human ad-hocery.

$\bullet$ The constant $k$ such that $p(n)=O(k^n)$, $q(n)=O(k^n)$, that happens to be the
largest root of the equation in $N$ obtained by equating to $0$ the leading coefficient in $n$ of $ope(n,N)$.
This follows easily from the so-called Poincar\'e lemma (see [vdP]).

$\bullet$ The constant $\beta$ such that $c(n)=O(\beta^n)$, analogously obtained from
$ope_1(n,N)$.

$\bullet$ The implied $\delta$ computed by
$$
\delta := \frac{2 \ln |k|  -\ln | \beta | }{ \ln |k| + GR + \ln P} \quad,
$$
that follows from the fact, equivalent to the prime number theorem, that $L(n)=O(e^n)$.

$\bullet$ The implied {\it irrationality measure}$\frac{\delta}{\delta -1}$  (if $\delta<1$ then it is negative, of course).

$\bullet$ The floating-point approximation for our constant $\alpha:= \lim_{n \rightarrow \infty} \frac{p(n)}{q(n)}$.

$\bullet$ If Maple succeeds to {\tt identify} that constant in terms of well-known ones, its conjectured
identification, otherwise, it just gives it back.

(Note, this part is only conjectural for now)

$\bullet$ The initial conditions for $p(n)$

$\bullet$ The initial conditions for $q(n)$

The output of {\tt RAsum} supplies everything needed for a verbose, human-readable, proof, done
by procedure {\tt RAsumV}.

{\bf Some Examples of RAsum with $\delta>1$}

Typing, in {\tt NesApery}

{\tt lprint( RAsum(binomial(n,k)*binomial(n+k,k)*3**k,k,n,N,10));}

immediately yields the output

{\tt
[(-n-2)*N**2+(14*n+21)*N-n-1, -n/(n+1)+N, [1, 1, 1, 1], [7+4*3**(1/2), 13.92820323], 1, 
[2*ln(7+4*3**(1/2))/(ln(7+4*3**(1/2))+1), 1.449629514], [2*ln(7+4*3**(1/2))/(ln(7
+4*3**(1/2))-1), 3.224053290], .2876820724517809274392190059938274315035, ln(4/3), [0, 2], [1, 7]] \quad ,
} 

proving that $\ln(4/3)$ is irrational with $\delta= 1.449629514$.

In fact, for any integer $a \geq 1$, the sum
$$
\sum_{k=0}^{n} {{n} \choose {k}} {{n+k} \choose {k}} a^k \quad ,
$$
yields the irrationality of $\ln((a+1)/a)$, as first proved by the  humans
Krishna Alladi and M.L. Robertson (see [AR]), with
$$
\delta=
{\frac { 4 \, \ln  \left( \sqrt{a} + \sqrt{a+1} \right) }{2 \, \ln  \left(  \sqrt{a} + \sqrt{a+1} \right) +1}} \quad .
$$
This was nicely extended by human Laurent Habsieger(see [H]).

As an example of one of the many disappointments, try:

{\tt lprint(RAsum(binomial(n,k)**2*binomial(n+k,k)*2**k,k,n,N,10));}

You would get a certain third order operator, a certain unidentified constant,   $0.02266573727755793921\dots$, and
the disappointing $\delta= 0.4934093797$.

\bigskip

{\bf Constant Terms and the Amazing Almkvist-Zeilberger Algorithm}

Recall that the ``Constant Term'' functional, $CT_x$,  is defined as follows
$$
CT_x [P(x)] \, := \, Coefficient \quad of \quad x^0 \quad in \quad  P(x) \quad ,
$$
where $P(x)$ is any Laurent polynomial in $x$.

Let's look again at the denominator sequence for $\ln  ((a+1)/a)$, namely
$$
q(n)=\sum_{k=0}^{n} {{n} \choose {k}} {{n+k} \choose {k}} a^k \quad .
$$

Since, of course (thanks to the binomial theorem)
$$
{{n+k} \choose {k}}= CT_x \left [ \frac{(1+x)^{n+k}}{x^k} \right ] \quad,
$$
we have
$$
q(n)=\sum_{k=0}^{n} {{n} \choose {k}}  CT_x \left [ \frac{(1+x)^{n+k}}{x^k} \right ]  a^k  \quad .
$$
Since $CT$ is a {\bf linear functional}, we can put it in front of the summation sign, getting
$$
q(n)=  CT_x \left [ \sum_{k=0}^{n} {{n} \choose {k}}  \frac{(1+x)^{n+k}}{x^k}  a^k  \right ]  \quad .
$$
$$
=  CT_x \left [ (1+x)^n \sum_{k=0}^{n} {{n} \choose {k}}  \left ( \frac{a(1+x)}{x} \right )^k   \right ]  \quad .
$$
$$
=  CT_x \left [ (1+x)^n  \left ( 1+ \frac{a(1+x)}{x} \right )^n   \right ]  \quad .
$$
$$
=  CT_x \left [ \left ( \frac{(1+x)(a+(a+1)x)}{x} \right )^n   \right ]  \quad .
$$

So let's forget about binomial coefficients sums (alias terminating hypergeometric series) and look for Ap\'ery miracles arising from
denominators $q(n)$ given by
$$
q(n)= CT_x \left [ P(x)^n \right ] \quad,
$$
for {\it any} Laurent polynomial $P(x)$ with integer coefficients.

Luckily we have the {\bf Almkvist-Zeilberger algorithm} ([AZ])  (also included in {\tt NesApery})
that, in a fraction of a second can find a linear recurrence
operator annihilating $q(n)$, that we know {\it a priori} are integers. Then we proceed as before,
define a sequence $p(n)$, with different initial conditions,
that unlike the $q(n)$ are no longer integers but rather
rational numbers, then define $\alpha :=\lim_{n \rightarrow \infty} \frac{p(n)}{q(n)}$, and look for
miracles. The analogous procedures in {\tt NesApery}
are 

$\bullet$ {\tt DelSeqCT}: for the empirical $\delta$ and possible identification.

$\bullet$ {\tt RAct}: The Almkvist-Zeilberger analog of {\tt RAsum} (that uses Zeilberger). See quite a few examples in:

{\tt http://www.math.rutgers.edu/\~{}zeilberg/tokhniot/oNesApery10} \quad .

$\bullet$ {\tt RActV}:  verbose, humanly-readable version of  {\tt RAct}. See quite a few examples in:

{\tt http://www.math.rutgers.edu/\~{}zeilberg/tokhniot/oNesApery11} \quad .

For example for {\it any} positive integers $a<b$

{\tt RAct((1+a*x)*(1+b*x)/x,x,n,N,10)}

yields a sequence of rational approximations for $\ln (b/a)$ with 
$$
\delta=
2\, 
\frac{2 \ln    ( \sqrt{a}+\sqrt{b} )  -\ln  ( b-a )}
{ 2 \ln  ( \sqrt{a}+\sqrt{b}  ) +1}
$$

that for many $a$'s and $b$'s is larger than $1$.

{\bf An Ap\'ery Miracle with a Third-Order Recurrence}

The family of constant term sequences
$$
K(a,b,c)(n):= CT_x \left [  \left ( \frac{(1+ax)(1+bx)(1+cx)}{x^2} \right )^n  \right ]
$$

yields Ap\'ery miracles with $\delta>1$ for many choices of integers $a,b,c$, even though the
recurrence satisfied by the $q(n)$ (and of course $p(n)$) is {\bf third order}.
Alas we have no clue whether the defined constants, $\alpha:=\lim_{n \rightarrow \infty} \frac{p(n)}{q(n)}$,
are `interesting'.

{\bf Integral Representations}

The drawback of {\it both} binomial coefficients sums (alias terminating hypergeometric series) and 
`constant terms' is that
we have to first {\it guess} the constant, and then prove our guess. 
There is another way to get a sequence of rational approximations to constants, where
the constant is known {\it a priori}.

Consider (see Wadim Zudilin's beautiful paper [Zu]) the sequence
$$
I(n) := \int_0^1 \, \left ( \frac{x(1-x)}{1-(1-a)x)} \right )^n \frac{1}{1-(1-a)x} \, dx \quad ,
$$
where $a>0, a \neq 1$ is rational.

Using the Almkvist-Zeilberger algorithm once again, it is immediate to get a 
second-order linear recurrence
equation with polynomial coefficients 
(these polynomial coefficients, in turn, have integer coefficients).
Since $I(0)$ and $I(1)$ are obviously linear combinations with rational coefficients of $1$ and $\ln a$,
it follows (using the recurrence and induction)
that so does $I(n)$ for all $n$. Since $I(n)$ is small for large $n$, pretending that
it is $0$ yields better and better rational approximations for $\ln a$ and also imply
rigorous irrationality measures with $\delta$ often larger that $1$. By considering more general integrals and clever penny-pinching,
E. Rukhazde [Ru]  got even better irrationality measures, as beautifully narrated in [Zu].

These are handled by procedure {\tt RAint} for the terse version and {\tt RAintV} for the verbose version.

For a few samples see:

{\tt http://www.math.rutgers.edu/\~{}zeilberg/tokhniot/oNesApery12} (terse version) \quad ,

{\tt http://www.math.rutgers.edu/\~{}zeilberg/tokhniot/oNesApery13} (verbose version) \quad .

{\bf Irrationality Measures for arctan(1/k)}

For $arctan(1/(2k+1))$ we use the integral
$$
2^{3n}(2k+1)^{2n}
\int_0^1 \,  \frac{x^{2n}(1-x)^{2n}}{((2k+1)^2+x^2)^{2n+1}} \, dx \quad ,
$$
and RAint gave us
$$
\delta= 
2\,{\frac {\ln  \left( \, 16\,{k}^{2}+16\,k+ 6+ 8\,\sqrt {2}  \, k \,\sqrt {2\,{k}^{2}+2\,k+1}+ \, 4\,\sqrt {2}\sqrt {2\,{k}^{2}+2\,k+1} \, \right) 
-\ln 2 }{\ln  \left( \, 16\,{k}^{2}+16\,k+\, 6+ \, 8\,\sqrt {2}  \, k \, \sqrt {2\,{k}^{2}+2\,k+1}+4\,\sqrt {2}\sqrt {2\,{k}^{2}+2\,k+1}  \, \right) +2
}} \quad ,
$$
that is larger than $1$ for $k>0$. When $k=0$, we get $\delta=0.79119792..$, so we can't deduce the irrationality of $\pi$ from this.

For $arctan(1/(2k))$ we use a similar integral (with different normalization, see the package), getting
$$
\delta=
2\,{\frac {\ln ( \, 4\,k+2\,\sqrt {4\,{k}^{2}+1}  \, )\, - \, \ln 2 }{\ln  (  \, 4\,k+2\,\sqrt {4\,{k}^{2}+1} \, ) +1}} \quad ,
$$

that is larger than $1$ for $k \geq 2$.

{\bf Conclusion}

This article should be considered as yet-another {\bf case study} in teaching computers human tricks and methods,
and let them do {\bf symbol-crunching} rather than mere {\bf number-crunching}. We have only
scratched the surface, and we hope that other, more persistent, and possibly smarter people (and computers!),
will use the package {\tt NesApery} to discover further new and exciting Ap\'ery miracles.
Better still, discover new kinds of miracles, and then using the present {\bf methodology},
{\bf teach} computers how to find these {\it new kind of miracles}. {\it Amen ken yehi ratson}.

{\bf References}

[AR] K. Alladi and M. L. Robertson, {\it Legendre polynomials and irrationality},
J. Reine Angew. Math. {\bf 318}(1980), 137-155.

[AZ] Gert Almkvist and Doron Zeilberger,
{\it The method of differentiating under the
integral sign}, J. Symbolic Computation {\bf 10}(1990), 571-591; \hfill\break
{\tt http://www.math.rutgers.edu/\~{}zeilberg/mamarimY/duis.pdf} 

[Ba] A. Baker, {\it Rational approximations to certain algebraic numbers}, Proc. London Math. Soc. {\bf 14}(1964), 385-398.

[Be] Michael A. Bennett, {\it Effective measures of irrationality for certain algebraic numbers},
Journal of the Australian Mathematical Society (Series A),{\bf 62}(1997), 329-344;  \hfill\break
{\tt www.math.ubc.ca/\~{}bennett/paper7.pdf} \quad .

[H] Laurent Habsieger, {\it Linear recurrent sequences and irrationality measures},
J. Number Theory {\bf 37}(1991), 133-145; \hfill\break
{\tt https://www.ima.umn.edu/preprints/Jan89Dec89/489.pdf} \quad .

[N] Ian Niven, ``{\it Irrational Numbers}'',  The Carus Mathematical Monographs {\bf \#11},
The Mathematical Association of America, 1956. Third printing, September, 1967.

[vdP] Alfred  van der Poorten, {\it A proof that Euler  missed ... Ap\'ery's proof of the irrationality of $\zeta(3)$,
An informal report}, The Mathematical Intelligencer {\bf 1}(1979), 195-203; \hfill\break
{\tt http://www.ega-math.narod.ru/Apery1.htm} \quad .

[Ru] E. A. Rukhadze, {\it A lower bound for the approximation of $\ln 2$ by rational numbers},
Vestnik Moskov. Univ. Ser. I Mat. Mekh. [Moscow Univ. Math. Bull.] {\bf 42}(1987), 25-29, (In Russian).

[Ze] Doron Zeilberger, {\it The method of creative telescoping},  J. Symbolic Computation {\bf 11}(1991), 195-204; \hfill\break
{\tt  http://www.math.rutgers.edu/\~{}zeilberg/mamarimY/ct.pdf} \quad  .

[Zu] Wadim Zudilin, {\it An essay on irrationality measures of pi and other logarithms}, 
 Chebyshevskii Sbornik {\bf 5}(2004), Tula State Pedagogical University, 49-65 (Russian); \hfill\break
English version: {\tt http://arxiv.org/abs/math/0404523} .
\bigskip
\bigskip
\hrule
\bigskip
\bigskip
Shalosh B. Ekhad, c/o D. Zeilberger, Department of Mathematics, Rutgers University (New Brunswick), Hill Center-Busch Campus, 110 Frelinghuysen
Rd., Piscataway, NJ 08854-8019, USA. 
\bigskip
Doron Zeilberger, Department of Mathematics, Rutgers University (New Brunswick), Hill Center-Busch Campus, 110 Frelinghuysen
Rd., Piscataway, NJ 08854-8019, USA. \hfill\break
url: {\tt http://www.math.rutgers.edu/\~{}zeilberg/}   
\quad . \hfill\break
Email: {\tt zeilberg at math dot rutgers dot edu}   \quad .
\bigskip
\hrule
\bigskip
{\bf EXCLUSIVELY PUBLISHED IN THE PERSONAL JOURNAL OF SHALOSH B. EKHAD and DORON ZEILBERGER} 
{\tt http://www.math.rutgers.edu/\~{}zeilberg/pj.html} and {\tt arxiv.org} \quad .
\bigskip
\hrule
\bigskip
May 17, 2014
\end